\documentclass[journal, dvipsnames, fleqn]{IEEEtran}

\usepackage{graphicx, xcolor, here}
\usepackage{cite}
\usepackage{amsmath}
\usepackage{amsfonts}
\usepackage{booktabs}
\usepackage{dsfont}
\usepackage{amssymb}
\usepackage{url}
\usepackage{textcomp, algorithm}
\usepackage{algpseudocode}
\usepackage{xcolor}
\usepackage{listings}

\newcommand{\e}{\vskip 2mm}

\usepackage{tikz}     
\usetikzlibrary{arrows,shapes,positioning,calc,backgrounds,decorations.markings,quotes,angles}

\title{A Unified Efficient Gradient-Based Heuristic For Box-Constrained Expectation-Related and Risk-Averse Stochastic Optimization Problems}
\author{Mazen Alamir\thanks{The author is with Univ. Grenoble Alpes, CNRS, Grenoble INP, GIPSA-lab, 38000 Grenoble, France. Email: mazen.alamir@grenoble-inp.fr. Website: https://www.mazenalamir.fr.}}
\date{\today}
\begin{document}

\maketitle

\begin{abstract}
This paper presents a new algorithm addressing the problem of stochastic optimization where the cost function depends on a vector of uncertain parameters with known statistics. The algorithm is parameterized so as to address various stochastic formulations spanning from Expectation-focused to Value-at-Risk (VaR) as well as Conditional-Value-at-Risk (CVaR)-focused formulations. The algorithm leverages a recently proposed gradient-based Search \& Accelerate algorithm which is originally dedicated to deterministic optimization problems. The approach is based on a sequence of \texttt{warm-started} solutions of instances of the problem. These solutions together with a samples of other solutions belonging to the convex hull of the first ones constitute the set of admissible candidates. Among this discrete set of candidates, the optimal solution is selected with regards to a sample-based approximation of the targeted criterion. The relevance of the algorithm and its efficiency are discussed and shown using a tailored illustrative example. 
\end{abstract}


\section{Introduction}
There is no single formulation of stochastic optimization under uncertainties. The formulation depends on the nature of the problem and the degree of its sensitivity to risk. Three main formulations might be considered that are briefly introduced in the next section.
\subsection{The stochastic optimization formulations}
Consider a parameterized cost function $f$ that depends on a vector of decision variables $x\in \mathbb X\subset \mathbb R^n$ and a vector of uncertain parameters $p\in \mathbb R^{n_p}$. Stochastic optimization \cite{RoysetJohannesO.2025RAtS} is about finding the \textit{best} decision variable with regards to some \textbf{measure}, denoted hereafter by $J(x)$. This measures describes a set of outcomes, as expressed by the possible values of $f(x,p)$, when $p$ spans a future unknown set $\mathcal R$ of realizations, namely:
\begin{equation}
\mathcal F_\mathcal{R}(x) := \Bigl\{f(x,p)\Bigr\}_{p\in \mathcal R} \label{defdecalF}
\end{equation}
The most widely considered problem is the one that is expectation-oriented. In this formulation one considers the asymptotic case where $\mathcal R:=\mathbb R^{n_p}$ with the measure being the expectation of the cost function's value:
\begin{equation}
J_E(x):=\mathbb E\bigl[f(x,\cdot)\bigr] := \int_{\mathcal{R}=\mathbb R^{n_p}}\pi(p)f(x,p)dp \label{expectation}
\end{equation}
where $\pi$ is the probability density function of the random variable $p$. This measure does not take into account the dispersion of the outcomes around the expected value which might be incompatible with risk-averse applications such as cancer treatment \cite{moicancer}, safety-involving engineering issues in aeronautics \cite{SERES20232013, Morales2021GradientBE} or some specific financial portfolio optimization to cite but few examples. 
\e 
This is the reason why many works (see \cite{GaivoronskiAlexei2005Vipo, Shapiro} and the references therein) consider the measure to be the value-at-risk (\texttt{VaR}) at level $\alpha$ which corresponds to the $\alpha$-quantile of $\mathcal F_{\mathbb{R}^{n_p}}(x)$, namely:
\begin{equation}
J_V(x,\alpha):=\texttt{VaR}_\alpha(x) := \inf\bigl\{z:\ \texttt{Pr}[f(x,p)\le z]\ge \alpha\bigr\} \label{defdeVaR}
\end{equation}
in which the probability can be reformulated in a similar way to \eqref{expectation}:
\begin{equation}
\texttt{Pr}[f(x,p)\le z] := \int_{\mathcal R:=\mathbb R^{n_p}} \pi(p)\bigl[\mathds{1}_{f(x,p)\le z}\bigr]dp \label{}
\end{equation}
Here again, for highly risk-averse contexts, \texttt{VaR} measure might still be insufficiently \textit{careful} as it does not consider higher quantiles (beyond $\alpha$). Here is where the Conditional Value at Risk (\texttt{CVaR}) might be useful \cite{Larsen2002} as it is defined to be the expectation of the quantiles that are higher than $\alpha$, namely: 
\begin{equation}
J_C(x,\alpha):=\texttt{CVaR}_\alpha(x) := \mathbb E\bigl[L\ \vert\ L\ge \texttt{VaR}_\alpha(x)\bigr] \label{defdeCVaR}
\end{equation}
Shortly speaking, $\texttt{VaR}_\alpha$ expresses the $\alpha$-tail \textit{threshold} while $\texttt{CVaR}_\alpha$ expresses the $\alpha$-tail \textit{severity}. Obviously by definition $\texttt{CVaR}_\alpha\ge \texttt{VaR}_\alpha$.
\e 
From a computational perspective, $\texttt{CVaR}$ is appealing compared to \texttt{VaR}. The reason  is that, while \texttt{VaR} is generally non-smooth and difficult to optimize directly \cite{Larsen2002}, it can be shown that \texttt{CVaR} can be computed using\footnote{Here $(z)_+:=\max(z,0)$.}:
\begin{equation}
\texttt{CVaR}_\alpha(x)=\min_{\eta}\Bigl[\eta-\dfrac{1}{1-\alpha}\mathbb E\bigl[(f(x,\cdot)-\eta)_+\bigr]\Bigr] \label{alternativeformulation}
\end{equation}
\e The main algorithmic options  that address the three above mentioned stochastic optimization problems are briefly discussed and some representative references are provided in order to clearly justify and position the proposed algorithm. 
\subsection{Algorithms: a brief state of the art}
Let us first of all underline the fact that stochastic optimization corresponds to a rather huge literature and an exhaustive analysis of existing works go obviously beyond the scope of this contribution. Rather, a few representative references are provided and a rough clustering of the corresponding alternatives is proposed in order to precisely state and justify the current contribution. 
\e 
\subsubsection{The expectation-related algorithms}\label{sec-expectation}
The most intuitive although heavy to compute solution of \eqref{expectation} is to use a sample-based approximation of the cost function and to use any standard optimization algorithm such as BFGS, SQP, Newton, etc. in order to solve the resulting deterministic optimization problem. This leads to the so-called Sample-Average-Approximation approach \cite{Kleywegt2002SAA}. Namely, a high number $N$ of samples of the uncertain parameter is randomly drawn and the cost function defined by the average $1/N\sum_{k=1}^Nf(x,p^{(k)})$ is used. The drawback of this approach is that the expensive average is involved in each single iteration (with a new iterate $x$) which might be too heavy when large values of $N$ are needed to achieve a decent approximation. 
\e This triggered the interest in the \textbf{Stochastic Gradient} option in which, a single new random value of $p$ is triggered at each iteration of a gradient method \cite{Nemirovski2009Robust}. The rationale is that this process produces ultimately a set of iterates with an average that is close to the optimum following the Polyak-Ruppert averaging principle \cite{PolyakJuditsky1992}. It is not clear however whether the stochastic gradient scheme can accommodate for accelerated gradient-related tricks that build on the value of the previous iteration. 
\e 
\subsubsection{\texttt{VaR}-related algorithms}
The application of the sample-based average approach to the \texttt{VaR} optimization problem leads to a cost function that is typically non smooth and can be locally flat or discontinuous in its derivatives which induces troubles in the standard optimization algorithms. This is the reason why the major works on this problem adopt the paradigm of \textit{chance constrained optimization} \cite{GENG2019341, KUCUKYAVUZ2022100030, Pagnoncelli2009SAAChance, Luedtke2008Probability, alamo2009randomized} in which solutions that maximize $z$ are searched that meet some tractable probabilistic constraints such as: 
\begin{equation}
\sum_{k=1}^{N(\eta,\varepsilon,m)}\mathds{1}_{\bigl[f(x,p^{(k)})\ge  z\bigr]}\le m \label{probab}
\end{equation}
where $\eta$ and $\epsilon$ stands for the so-called precision and confidence parameters\footnote{These are typically small values. $\eta$ stands for the probability of the inequality being violated while $1-\varepsilon$ represents the confidence one might have in the previous . Typically, the pair $(\eta,m)$ induces the necessary condition $\eta N\ge m$ \cite{alamo2009randomized}.}, $m$ is a chosen integer while $N(\eta,\varepsilon,m)$ represents the number of samples $p^{(k)}$ of the uncertain parameters that has to be randomly drawn to be involved in the inequality. Notice however that the associated cost function is not smooth and hence depending upon structure, formulations therefore introduce binary variables, decomposition methods, relaxations, or other representations.
\e In order to avoid non smoothness in the quantile estimation that is related to the final number of samples, the approach in \cite{Hong2009Quantile} uses a finite number of samples in order to build an estimation of the sensitivity of the quantile to the decision variable since the original rigorous (non approximated) definition of \texttt{VaR} might be differentiable. This being said, the quite heavy process of estimating the sensitivity of the $\alpha$-quantile to $x$ is to be performed at each iteration.
\e For all the previous reasons, the \texttt{VaR}-related problems seem to be the least tractable of the three formulations. Indeed  the one associated to the \texttt{CVaR} is much more tractable. 
\e 
\subsubsection{\texttt{CVaR}-related algorithms}
Notice that by virtue of \eqref{alternativeformulation}, the \texttt{CVaR} cost takes the form of an expectation \cite{RockafellarUryasev2002CVaR}. Consequently, the same approaches described earlier in Section \ref{sec-expectation} namely, the sample-averaging approach \cite{Tamar2015CVaR} and the sensitivity based approach \cite{HongLiu2009CVaR} are eligible. Moreover, by jointly optimizing $x$ and $\eta$ involved in \eqref{defdeCVaR}, the corresponding \texttt{VaR} is precisely the value of the optimal $\eta$. 
\e 
\noindent \underline{\texttt{Discussion}}
\e 
From the above quick overview, it comes out that the available alternatives adopt one of two approaches, namely:\\
\begin{itemize}
\item  \textbf{Either} a \textbf{sample-based approximation} of some expectation that is then used to defined the cost function. This is the case when the expectation formulation is addressed  or when the \texttt{CVaR}-related formulation is addressed thanks to the Rockafellar–Uryasev representation \eqref{alternativeformulation}. This option leads to rather heavy-to-compute cost functions when high values of the number of samples $N$ is adopted in order to get a relevant approximation of the expectation being involved for each current iterate $x$ during the optimization process.\\
\item \textbf{or} a \textbf{stochastic gradient} approach is adopted in which each step uses a gradient direction that is associated to a randomly sampled value of the uncertain parameter. This induces a much computationally cheap iterations during which a number of iterates is collected and some of the \textit{last} ones are aggregated to get the final solution. The relevance of the resulting aggregated solution is conditioned by the number of iterations being \textit{sufficiently} high. Moreover, the fact that the involved uncertain parameter is modified at each iteration makes risky the use of \textit{acceleration}-oriented versions of the gradient descent such as the one adopted in \cite{nesterovfg2004,alamir2026nonlinearmodelpredictivecontrol} and the related references. This might lead to the use of the notoriously slow simple raw gradient method. This with the need for a high number of iterations are the main drawbacks of this approach. \\
\end{itemize}
Combining the two sets of issues associated to the two approaches, it is likely that the resulting algorithms might be incompatible with real-time control/estimation scheme arising in Stochastic Nonlinear Model Predictive control (\textbf{SNMPC}) \cite{Mesbah2016} or Stochastic Nonlinear Moving Horizon Estimation) (\textbf{SNMHE}) \cite{Aghapour2020RiskMHE, muntwiler2025mheparametricuncertainty}, at least for a class of systems requiring high updating rates. This justifies the relevance of the algorithm proposed in the present contribution. 
\e More precisely, in the proposed algorithm, while a sample-based approximation is ultimately used, it is not involved in the inner optimization loop, rather, a set of candidate values is created by first efficiently solve (through a warm-start mechanism) a set of deterministic optimization problems. The so-collected solutions are then used to generate some extra ones from their convex hull and it is only then that the sample-based approximation is used in order to choose the optimal one among all of the candidate solutions. This process can then be repeated to progressively improve the quality of the solution. This enable to reduce the computational burden while leveraging the efficient acceleration-empowered gradient-based iterations proposed in \cite{alamir2026nonlinearmodelpredictivecontrol}. 
\e 
The remainder of this paper is organized as follows: First Section \ref{sec-algo} progressively explains the proposed algorithm. The stochastic optimization problem that is used in the numerical investigation of the proposed algorithm is described in Section \ref{sec-example} together with all the settings used in the numerical experiments. The results are shown in Section \ref{sec-results}. Finally section \ref{sec-conc} concludes the paper and give some hints for further investigation.
\section{Description of the proposed algorithm}\label{sec-algo}
Before the rigorous code of the algorithm is provided, its main steps are first explained hereafter: \\
\begin{enumerate}
\item First the algorithm generates a discrete set of $\texttt{nDesign}$ candidate values, say $x^{(i)}_\star$, $i=1,\dots,\texttt{nDesign}$. This is done by warm-start solving a sequence of problems using the very efficient acceleration-endowed gradient-based algorithm, recently proposed in \cite{alamir2026nonlinearmodelpredictivecontrol}. This set of candidate values are referred to hereafter by $\mathbb X_0^\star\subset \mathbb R^n$:
\begin{equation}
\mathbb X_0^\star = \Bigl\{x_\star^{(i)}\mid i=1,\dots, \texttt{nDesign} \Bigr\} \label{defdeXmathbbstar}
\end{equation}
\e More precisely, a set of $\texttt{nDesign}$ instances of the uncertain parameters (preferably containing the nominal one, if any, in the first position) is randomly sampled. Then the list is iteratively reordered starting from the second one so that each element is the closest to the previous one among all the still not-yet-ordered elements. The resulting set is denoted as follows: 
\begin{subequations}
\begin{align}
&\mathcal P_0:=\Bigl\{p^{[1]}, p^{[2]}\dots, p^{[\texttt{nDesign}]}\Bigr\} \quad\text{$p^{[1]}$: nominal}\label{defdecalP} \\
&p^{[i+1]} = \text{arg}\min_{p\in \mathcal{P}_0-\{p^{[j]}\mid j\ge i+1\}}\|p-p^{[i]}\| \label{reorderP0}
\end{align}    
\end{subequations}
Once this ordering is done, the set of optimization problems associated to $p^{[i]}$ are successively solved with warm starting being used for all $i\ge 2$:
\begin{subequations}
\begin{align*}
x_\star^{(1)}&\leftarrow \texttt{solve}(f, x_\text{initial}, p^{[1]}, n_1)\label{solve1}\\
x_\star^{(i+1)}&\leftarrow\texttt{solve}(f, x_\star^{(i)}, p^{[i]}, n_2)
\end{align*}    
\end{subequations}
where the notation in the r.h.s of the above notation means:
\begin{center}
\tikz{
\node at(0,0) (O){\texttt{solve}$(f,x_0,p,n)$};

\node[anchor=west] at($(O.east)+(0.2,0)$){
\begin{minipage}{0.25\textwidth}
Use the deterministic algorithm \texttt{solve} to minimize (in $x$) the cost function $f(x,p)$ using $x_0$ as initial guess and allowing a maximum number of $n$ iterations
\end{minipage}
};
}
\end{center}
Notice that the ordering enables to increase the relevance of the warm starting mechanism by virtue of fact that each two successive values are supposed to be close to each other. This justifies the use of $n_2<n_1$ as warm starting enables less number of iterations to be required. 
\e In this contribution, the algorithm used inside the \texttt{solve} operator is the one recently proposed in \cite{alamir2026nonlinearmodelpredictivecontrol} which implements an efficient gradient-based algorithm that combines line search and acceleration mechanism. 
\\ 
\item The \textit{final} set $\mathbb X^\star$ of candidate solutions is obtained by adding to $\mathbb X_0^\star$, a set of randomly generated values that belong to the convex hull of $\mathbb X_0^\star$, namely: 
\begin{equation}
\mathbb X^\star := \mathbb X_0^\star \cup \mathbb X_{cv}^\star(\mathbb X^\star_0)\label{defdemathbbXstar}
\end{equation}
where 
\begin{subequations}\label{cvhull}
\begin{align}
X_{cv}^\star(\mathbb X^\star_0) &:= \bigcup_{k=1}^{n_\text{cv}}\bigl\{\lambda_k\xi_1^{(k)}+(1-\lambda_k)\xi_2^{(k)}\bigr\} \label{defdeXcv} \\
&=: \bigcup_{k=1}^{n_\text{cv}}\Bigl\{x_\star^{(k)}=:x_\star^{\texttt{nDesign}+k}\Bigr\}\label{lanotationstar+}
\end{align}
\end{subequations}
where for each $k$, $\lambda_k$, $(\xi_1^{(k)}, \xi_2^{(k)})$ are randomly generated respectively in $(0,1)$ and $\mathbb X_0^\star\times \mathbb X_0^\star$ such that $\xi_1^{(k)}$ and $\xi_2^\star$ are not identical. \\
\item Denoting by $n_c=\texttt{nDesign}+n_\text{cv}$ the number of elements in the resulting set $\mathbb X^\star$, we get a set of candidate values of the decision variable that, by virtue of the notation introduced in \eqref{lanotationstar+}, can be written as follows:
\begin{equation}
\mathbb X^\star = \Bigl\{x_\star^{(k)}\Bigr\}_{k=1}^{n_c} \label{defdeXstarnc}
\end{equation}
\item Given the set $\mathbb X^\star$ of candidate values, one can compute their associated \textit{merit} given a sample-based representation of any of the cost functions mentioned above, namely, that related to \texttt{Expectation}, \texttt{VaR} or \texttt{CVaR}. In order to do this, a number $N$ of samples of the uncertain parameters, say:
\begin{equation}
\mathbb P:=\Bigl\{p^{(i)}\Bigr\}_{i=1}^N \label{defdemathbbP}
\end{equation}

is randomly generated so that the following deterministic costs can be defined:
\begin{subequations}\label{tildaJ}
\begin{flalign}
\tilde J_E(x)&:=\dfrac{1}{N} \sum_{i=1}^Nf(x,p^{(i)})\label{lestildaJE}\\
\tilde J_V(x,\alpha)&:=\texttt{percentile}\Bigl(f(x,\mathbb P), \alpha\Bigr)\label{lestildaJV}\\
\tilde J_C(x,\alpha)&:=\texttt{Avg}\Bigl[v\in f(x,\mathbb P)\mid v\ge \tilde J_V(x,\alpha)\Bigr]\nonumber
\end{flalign}
\end{subequations}
Notice that it is possible at this stage to use a large number $N$ of samples as the approximated maps above are not involved in each iteration of an optimization process. Moreover, there is no need that $\tilde J_V(x,\alpha)$ be appropriate for optimization as the latter is done via enumeration over the discrete set $\mathbb X^\star$ of candidate values.
\e 
\item At this stage, for a given choice of the cost function among $\tilde J_E$, $\tilde J_V$ and $\tilde J_C$, the process can be terminated by taking the element of $\mathbb X^\star$ that minimizes the chosen cost, namely: 
\begin{equation}
x^\star_\ell(\mathbb X^\star) := \text{arg}\min_{x\in \mathbb X^\star} \tilde J_\ell(x)\mid  \ell\in \{E, V, C\}\label{defdeoptimalxXstar}
\end{equation} 
\item If there is still available computation time, it is possible to go further in the optimization by iterating over a several updates of the admissible set of candidates in order to attempt to find better solutions. 
\e In order to this, consider  that the previously discussed set $\mathbb X^*$ is simply an initial version of some set iteration, namely: 
\begin{equation}
\mathbb Z^{(0)} = \mathbb X^\star \label{defdeZ0}
\end{equation}
In order to define the next iterate and having the set of values corresponding to the elements of $\mathbb Z^{(0)}$, let us denote by $\mathcal B_m(\mathbb Z^{(0)})$ the set of first best $m$ candidates in terms of the cost function values. The next set to be examined, namely $\mathbb Z^{(1)}$ is defined by:
\begin{equation}
\mathbb Z^{(1)} := \mathbb X_\text{cv}^\star\bigl(\mathcal B_m(\mathbb Z^{(0)})\bigr)\in [\mathbb R^n]^{n_\text{cv}} \label{defdeZ1}
\end{equation}
It is then possible to use the iteration \eqref{defdeZ1} repeatedly:
\begin{equation}
 \mathbb Z^{(i+1)} := \mathbb X_\text{cv}^\star\bigl(\mathcal B_m(\mathbb Z^{(i)})\bigr)\in [\mathbb R^n]^{n_\text{cv}} \label{defdeZi}
\end{equation}
\begin{table*}
    \centering
    \caption{List of parameters involved in the algorithm}
    \begin{tabular}{|l|l|}
    \hline
        \textbf{Parameter} & \textbf{Description} \\
        \hline\hline
        $\texttt{nDesign}$ & Number of initial warm-start solved deterministic optimization problems. \\
        $n_\text{cv}$ & The number of selected candidate from the convex hull. \\
        $N$ & Number of samples for the approximation of the statistical functions. \\
        $m$ & Number of best selected candidate for the updating of $\mathbb Z^{(i)}$. \\
        \texttt{maxIterCV} & Maximum number of updating according to \eqref{defdeZi}. \\
        $n_1$ & Maximum number of iteration for \texttt{solve} for the first problem.\\
        $n_2$ & Maximum number of iteration for \texttt{solve} for the other problems.\\
        \hline 
    \end{tabular}
    \label{tabparam}
\end{table*}
hence enabling to examine potentially better population of candidates leading at each iteration to associate the best choice which, following \eqref{defdeoptimalxXstar}, is denoted by $x_\ell^\star(\mathbb Z^{(i)})$. \\
\item The set of best candidates that would be collected after the update is performed $\texttt{maxIterateCV}$ times is denoted by:
\begin{equation}
\mathcal X^\star :=  \Bigl\{x_\ell^\star(\mathbb Z^{(i)})\Bigr\}_{i=0}^\texttt{maxIterateCV}\label{defdeXcalstar}
\end{equation}
\item Hence the best solution finally delivered by the algorithm is defined by:
\begin{equation}
x^\text{opt}:=\text{arg}\min_{x\in \mathcal X^\star}\Bigl[\tilde J_\ell(x)\Bigr]\label{lastouf}
\end{equation}
where the optimization is simply done by reordering the already computed sub-optimal values $x_\ell^\star(\mathbb Z^{(i)})$, $i=1,\dots, \texttt{maxIterateCV}$. 
\end{enumerate}
The previous steps involve a set of parameters that can be summarized in Table \ref{tabparam} for the sake of conciseness. 
\e 
All the above described steps are summarized in Algorithm \ref{algo}. Moreover, Table \ref{tabparam} summarizes the signification of the main parameters used in the algorithm. 
\e 
\begin{algorithm}
\footnotesize
\caption{Gradient-based stochastic optimization}\label{algo}
\begin{algorithmic}[1]  
\Statex \textbf{Input parameters} (with suggested default values): 
\Statex $\texttt{nDesign}=20$, $n_\text{cv}=10$, $N=1000$, $m=10$, \texttt{maxIterCV}=5, $n_1=50$, $n_2=10$. \texttt{xinitial}$\in \mathbb R^n$. $\ell\in \{E,V,C\}$, $\alpha=0.95$.
\Statex \textbf{Initialization}: 
    \Statex $\mathcal P_0:=\bigl\{p^{[i]}\mid i=1,\dots, \texttt{nDesign}\bigr\}\leftarrow$ Generate $\texttt{nDesign}$ samples of parameters.
    \Statex $\mathcal P_0\leftarrow$ Reorder the elements \Comment{according to \eqref{reorderP0}.}
    \Statex $\mathbb X_0^\star \leftarrow  \emptyset$, 
    $n_c:=\texttt{nDesign}+n_\text{cv}$.
    \Statex \hrulefill
    \e 
    \Statex {\color{Gray} ----- Compute $\mathbb X_0^\star$ -----}
    \e 
    \State $x_\star^{(1)}, f_\star^{(1)}\leftarrow \texttt{solve}(f, \texttt{xinit}, p^{[1]},n_1)$ \Comment{$f_\star^{(1)}$ optimal value}
    \For{$i=2,\dots,\texttt{nDesign}$}
    \State $x_\star^{(i)}, f_\star^{(i)}\leftarrow\texttt{solve}(f, x_\star^{(i-1)}, p^{[i]}, n_2)$ \Comment{$f_\star^{(i)}$ optimal value}
    \EndFor
    \State $\mathbb X_0^\star\leftarrow \bigl\{x_\star^{(i)}\bigr\}_{i=1}^{\texttt{nDesign}}$
    \e 
    \Statex {\color{Gray} ----- Extend $X_0^\star$ by random sampling inside its convex hull -----}
    \e 
    \State Sample $\mathbb X_\text{cv}(\mathbb X_0^\star)=:\bigl\{x_\star^{(i)}\bigr\}_{i=\texttt{nDesign}+1}^{n_c}$ \Comment{According to \eqref{cvhull}}
    \e 
    \Statex {\color{Gray} ----- Find by enumeration the optimal solution and cost -----}
    \e 
    \State $f_\star^{(i)}=\tilde J_\ell(x^{(i)}, \alpha)$, $i\in \{\texttt{nDesign}+1, \dots, n_c\}$
    \State $i^\text{opt}\leftarrow$ arg$\displaystyle{\min_{i=1,\dots,n_c}} f_\star^{(i)}$
    \State $\texttt{sol}^{(0)}\leftarrow (x_\star^{i^\text{opt}}, f_\star^{i^\text{opt}})$
    \e 
    \Statex {\color{Gray} ----- Try to improve by selecting inside convex hull of best solutions -----}
    \e 
    \State $Z^{(0)}\leftarrow \mathbb X_0^\star\bigcup \mathbb X_\text{cv}(\mathbb X_0^\star)$
    \For{$k=1,\dots,\texttt{maxIterCV}$}
    \State $\mathbb Z^{(k)}\leftarrow \mathbb X_\text{cv}\bigl(\mathcal B_m(\mathbb Z^{(k-1)})\bigr)=:\bigl\{x^{(i)}\bigr\}_{i=1}^{n_\text{cv}}$
    \For {$i=1,\dots, n_\text{cv}$}
    \State $f^{(i)}=\tilde J_\ell(x^{(i)}, \alpha)$
    \EndFor
    \State $i^\text{opt}\leftarrow$ arg$\displaystyle{\min_{i=1,\dots,n_\text{cv}}} f^{(i)}$
    \State $\texttt{sol}^{(k)}\leftarrow (x^{i^\text{opt}}, f^{i^\text{opt}})$
    \EndFor
    \e 
    \Statex {\color{Gray} ----- Determine the returned solution -----}
    \e 
    \State \textbf{Return} $\texttt{sol}^{(i)}$ corresponding to the lowest second attribute $f$. 
\end{algorithmic}
\end{algorithm}
In the next section, the example that serves in the numerical experiments is described. 
\section{Illustrative stochastic optimization problem}\label{sec-example}
Consider the following instantiation of the cost function: 
\begin{subequations}
\begin{align}
f(x,p)&:= \varphi(x) + \rho \Psi(p) e^{-\lambda \varphi(x)}\label{defdefex}\\
\varphi(x)&:= \|x-x_c\|^2 \quad;\quad x_c={\scriptsize \begin{bmatrix} 
0.0 \cr 0.25 \cr 0.5\cr 0.75\cr 1.25\cr 1.5
\end{bmatrix}} \label{defdevarphixex}\\
\Psi(p)&:= \|p-\mathbf 1\|^2 \label{defdevarPsixex}
\end{align}   
\end{subequations}
where $x\in \mathbb R^{6}$, $p\in \mathbb R^2$ and where $\mathbf 1$ stands for a vector of ones of the appropriate dimension. The values $\lambda=0.5$ and $\rho=50$ are used in the forthcoming experiments. The dispersion of $p$ is supposed to obey a normal distribution around $p_\text{nominal}=\mathbf 1$ with a standard deviation $\sigma=0.35$. 
\e The rationale in this definition is that the optimal value of $f(x, p)$ when $p=\mathbf 1$ is obviously $f(x^\star, \mathbf 1)=0$ with the optimal solution $x^\star=x_c$. As a matter of fact, the expected value of $p$ is precisely $\mathbf 1$ meaning that if the nominal most expected value of $p$ is used in a deterministic solution, the optimal solution would appear to be $x^\star=x_c$.
\e 
When the dispersion of $p$ is to be accounted for however, some values that are sufficiently far from $\mathbf 1$ would lead to a non vanishing second term in \eqref{defdefex} and this extra term would be even higher when the nominal optimal value $x=x_c$ is used. Consequently, a dispersion-aware optimal $x$ would be pushed away from $x=x_c$. This feature is obviously enforced with increasing $\rho$ and decreasing $\lambda>0$. 
\e 
In the following results, the solutions obtained by the proposed algorithm for the different criteria (\texttt{Expectation}, \texttt{VaR} and \texttt{CVaR}) will be compared to the solution based on the deterministic solution that is based on the nominal value $p_\text{nom}=\mathbf 1$. 
\e The set of \texttt{nDesign} parameters $\mathcal P_0$ used in the design [see \eqref{defdecalP}] and the set $\mathbb P$ used in the internal evaluation of the approximated version of the cost functions [see \eqref{defdemathbbP}], namely $\tilde J_E$, $\tilde J_V$ and $\tilde J_C$ are obtained by first generating $N$ samples of the parameters leading to $\mathbb P$ of which the first \texttt{nDesign} values are taken to define $\mathcal P_0\subset \mathbb P$. Recall that both sets are used in the computation of the \textit{sub-optimal} solution $x^\text{opt}$ defined by \eqref{lastouf}. 
\e Once a suboptimal solution is found, its performance is evaluated using a much larger \textit{test} set of instances of the uncertain parameters, denoted hereafter by $\mathbb P_\text{test}$ which is of cardinality $N_\text{test}\gg N$. 
\e 
As a matter of fact, since randomness lies in the process, the whole operations (generating $\mathbb P$, $\mathbb P_\text{test}$ and running the optimization algorithm and the test-related performance) is repeated (100 times) and statistics of the performance on the test set are drawn using percentiles of the achieved costs over all the test sets $\mathbb P_\text{test}$ encountered during these repetitions. These are the percentiles of the achieved costs and computation times which are reported on all the forthcoming Tables. 
\e 
In all the experiments, the number of instances included in the test set $\mathbb P_\text{test}$ is taken to be equal to $$\texttt{card}(\mathbb P_\text{test})=100000$$ while the number $N$ of values in the optimization set $\mathbb P$ and hence used in the optimization is taken such that $$\texttt{card}(\mathbb P)=N\in \{100, 1000\}$$ in order to show the impact of this choice on the quality of the resulting solutions as well as on the corresponding computation time. 
\section{Results}\label{sec-results}
Table \ref{resultsN50withoutN1000} shows the percentiles of the following quantities induced during the solution of the 100 random instances of problem (and its test sets $\mathbb P_\text{test}$) for each of the cost function's choice $\in \{\textbf{Expectation}, \textbf{VaR}, \textbf{CVaR}\}$:\\
\begin{description}
    \item[\textbf{Test: nominal}] \ \e This column gives the percentiles of the cost function values on the test dataset $\mathbb P_\text{test}$ for the solution obtained using the nominal value $p=\mathbf{1}$, or equivalently when the standard deviation $\sigma=0$ is used in the generation of $\mathbb P$. This columns serves as a basis for the appreciation of the uncertainty-aware solution provided by the algorithm. Notice that the parametric dispersion contained in the test dataset $\mathbb P_\text{test}$ is equal to $\sigma=0.35$.\e 
    \item[\textbf{Test: stochastic}] \ \e This column gives the percentiles of the cost function values on the test dataset $\mathbb P_\text{test}$ for the solutions obtained with a design set $\mathbb P$ involving parameters dispersion using the \textit{correct} standard deviation $\sigma=0.35$. \e 

    \item[\textbf{Train: stochastic}] \ \e This columns gives the percentile of the optimal values obtained during the optimization and computed using the comparatively small cardinality set $\mathbb P$. It represents the dispersion of the optimal values as predicted by the algorithm based on its small sampled working set of instances $\mathbb P$. \e 
    \item[\textbf{cpu (sec)}] \ \e This columns gives the percentile of the computation times needed by the algorithm to solve the 100 instances of the optimization problem for each of the cost function being targeted. 
\end{description}
\e 
While Table \ref{resultsN50withoutN1000} shows the results when the number of samples $N=1000$ is used to build $\mathbb P$, Table \ref{resultsN50withoutN100} shows the same results when $N=100$ samples is used in $\mathbb P$. Notice how the computation time is significantly reduced at the price of higher cost function for the very high quantiles of the results. Only the last two quantiles are impacted by the difference. 
\e 
While in the previously described tables \ref{resultsN50withoutN1000} and \ref{resultsN50withoutN100}, the number \texttt{nDesign} of deterministic optimization problem are first solved was equal to \texttt{nDesign}=50, Table \ref{resultsN20withoutN100} shows the results when \texttt{nDesign}=20 and $N=100$. Here again, the computation time is reduced approaching $0.04$ sec at the price a deterioration in the higher percentiles of the achieved costs. 
\e Notice however that the achieved statistics are still decently comparable showing the ability to monitor the trade-of between the risk mitigation level and the computation time which, as it is mentioned in the introduction, might be crucial and even unavoidable in the case of use of the proposed algorithm in a stochastic MPC framework where the computation time is a key parameter in the real-time implementability of these schemes.
\e Finally, in order to examine the benefit from using the extra-candidate values sampled inside the convex hull of the primary set of candidates, Table \ref{resultscvcompar} shows comparison of the statistics with and without this specific step in the algorithm. The results suggest that, despite a slight improvement of the quality (for intermediate quantiles for \textbf{Expectation}-and \textbf{CVaR}-related costs, it seems that the benefit is rather weak, as far as the problem under consideration is concerned and given the significant induced increase in the computation times. 
\begin{table*}[]
\begin{center}
\begin{tabular}{lcccr}
\toprule
 q\% | \textbf{Expectation} & Test: nominal & Test: stochastic & cpu (sec) & Train: stochastic \\
\midrule
0 & 12.16 & 5.61 & 0.50 & 5.49 \\
25 & 12.23 & 5.62 & 0.52 & 5.57 \\
50 & 12.25 & 5.63 & 0.52 & 5.63 \\
80 & 12.28 & 5.64 & 0.52 & 5.68 \\
90 & 12.30 & 5.64 & 0.53 & 5.70 \\
95 & 12.32 & 5.64 & 0.53 & 5.72 \\
99 & 12.33 & 5.64 & 0.55 & 5.80 \\
100 & 12.34 & 5.79 & 0.56 & 5.82 \\
\bottomrule
\end{tabular}
\e 
\begin{tabular}{lcccr}
\toprule
 q\% | \textbf{VaR} & Test: nominal & Test: stochastic & cpu (sec) & Train: stochastic \\
\midrule
0 & 36.30 & 7.80 & 0.51 & 7.62 \\
25 & 36.57 & 7.82 & 0.52 & 7.76 \\
50 & 36.71 & 7.83 & 0.52 & 7.82 \\
80 & 36.84 & 7.84 & 0.53 & 7.88 \\
90 & 36.94 & 7.87 & 0.53 & 7.91 \\
95 & 37.00 & 7.89 & 0.53 & 7.97 \\
99 & 37.10 & 8.06 & 0.54 & 8.04 \\
100 & 37.17 & 9.14 & 0.54 & 8.75 \\
\bottomrule
\end{tabular}
\e 
\begin{tabular}{lcccr}
\toprule
q\% | \textbf{CVaR} & Test: nominal & Test: stochastic & cpu (sec) & Train: stochastic \\
\midrule
0 & 48.35 & 8.38 & 0.52 & 8.10 \\
25 & 48.80 & 8.40 & 0.53 & 8.31 \\
50 & 48.94 & 8.41 & 0.53 & 8.41 \\
80 & 49.11 & 8.49 & 0.54 & 8.53 \\
90 & 49.26 & 8.53 & 0.54 & 8.59 \\
95 & 49.42 & 8.63 & 0.54 & 8.66 \\
99 & 49.60 & 8.93 & 0.55 & 8.89 \\
100 & 49.61 & 8.95 & 0.56 & 9.05 \\
\bottomrule
\end{tabular}

\end{center}    
\caption{Dispersion statistics (for the different costs) of the optimal solution over the 100 random trials. \\ \texttt{nDesign}=50, $\mathbf{N=1000}$}\label{resultsN50withoutN1000}
\end{table*}

\begin{table*}[]
\begin{center}
\begin{tabular}{lcccr}
\toprule
 q\% | \textbf{Expectation} & Test: nominal & Test: stochastic & cpu (sec) & Train: stochastic \\
\midrule
0 & 12.13 & 5.61 & 0.08 & 5.26 \\
25 & 12.23 & 5.63 & 0.09 & 5.48 \\
50 & 12.26 & 5.64 & 0.10 & 5.68 \\
80 & 12.28 & 5.64 & 0.10 & 5.81 \\
90 & 12.30 & 5.66 & 0.10 & 5.88 \\
95 & 12.31 & 5.67 & 0.10 & 5.99 \\
99 & 12.35 & 6.20 & 0.10 & 6.10 \\
100 & 12.35 & 6.36 & 0.10 & 6.61 \\
\bottomrule
\end{tabular}
\e 
\begin{tabular}{lcccr}
\toprule
 q\% | \textbf{VaR} & Test: nominal & Test: stochastic & cpu (sec) & Train: stochastic \\
\midrule
0 & 36.28 & 7.80 & 0.08 & 6.99 \\
25 & 36.57 & 7.82 & 0.09 & 7.56 \\
50 & 36.72 & 7.83 & 0.10 & 7.74 \\
80 & 36.85 & 7.87 & 0.10 & 8.00 \\
90 & 36.93 & 7.91 & 0.10 & 8.14 \\
95 & 37.00 & 7.98 & 0.10 & 8.26 \\
99 & 37.10 & 8.33 & 0.15 & 8.68 \\
100 & 37.11 & 36.67 & 0.19 & 32.56 \\
\bottomrule
\end{tabular}
\e 
\begin{tabular}{lcccr}
\toprule
q\% | \textbf{CVaR} & Test: nominal & Test: stochastic & cpu (sec) & Train: stochastic \\
\midrule
0 & 48.29 & 8.37 & 0.08 & 7.74 \\
25 & 48.74 & 8.40 & 0.09 & 8.17 \\
50 & 48.92 & 8.43 & 0.10 & 8.35 \\
80 & 49.13 & 8.54 & 0.10 & 8.68 \\
90 & 49.25 & 8.63 & 0.10 & 8.84 \\
95 & 49.35 & 8.76 & 0.10 & 9.01 \\
99 & 49.39 & 10.10 & 0.11 & 9.89 \\
100 & 49.40 & 38.66 & 0.12 & 40.29 \\
\bottomrule
\end{tabular}

\end{center}    
\caption{Dispersion statistics (for the different costs) of the optimal solution over the 100 random trials. \\ \texttt{nDesign}=50, $\mathbf{N=100}$}\label{resultsN50withoutN100}
\end{table*}

\begin{table*}[]
\begin{center}
\begin{tabular}{lcccr}
\toprule
 q\% | \textbf{Expectation} & Test: nominal & Test: stochastic & cpu (sec) & Train: stochastic \\
\midrule
0 & 12.15 & 5.61 & 0.04 & 5.16 \\
25 & 12.22 & 5.63 & 0.04 & 5.52 \\
50 & 12.25 & 5.63 & 0.04 & 5.65 \\
80 & 12.28 & 5.65 & 0.04 & 5.80 \\
90 & 12.31 & 5.67 & 0.04 & 5.84 \\
95 & 12.33 & 5.70 & 0.04 & 5.88 \\
99 & 12.33 & 5.76 & 0.04 & 6.01 \\
100 & 12.34 & 5.79 & 0.04 & 6.11 \\
\bottomrule
\end{tabular}
\e 
\begin{tabular}{lcccr}
\toprule
 q\% | \textbf{VaR} & Test: nominal & Test: stochastic & cpu (sec) & Train: stochastic \\
\midrule
0 & 36.36 & 7.80 & 0.03 & 7.06 \\
25 & 36.58 & 7.83 & 0.04 & 7.50 \\
50 & 36.72 & 7.86 & 0.04 & 7.70 \\
80 & 36.81 & 7.98 & 0.04 & 8.02 \\
90 & 36.97 & 8.09 & 0.04 & 8.22 \\
95 & 37.03 & 8.29 & 0.04 & 8.39 \\
99 & 37.07 & 9.59 & 0.04 & 8.57 \\
100 & 37.10 & 9.87 & 0.04 & 9.30 \\
\bottomrule
\end{tabular}
\e 
\begin{tabular}{lcccr}
\toprule
 q\% | \textbf{CVaR} & Test: nominal & Test: stochastic & cpu (sec) & Train: stochastic \\
\midrule
0 & 48.26 & 8.39 & 0.03 & 7.63 \\
25 & 48.75 & 8.43 & 0.04 & 8.29 \\
50 & 48.93 & 8.60 & 0.04 & 8.60 \\
80 & 49.11 & 8.99 & 0.04 & 9.06 \\
90 & 49.26 & 9.32 & 0.04 & 9.32 \\
95 & 49.28 & 9.78 & 0.04 & 9.53 \\
99 & 49.46 & 10.66 & 0.04 & 11.53 \\
100 & 49.50 & 16.37 & 0.04 & 16.19 \\
\bottomrule
\end{tabular}

\end{center}    
\caption{Dispersion statistics (for the different costs) of the optimal solution over the 100 random trials. \\  \texttt{nDesign}=20, $\mathbf{N=100}$}\label{resultsN20withoutN100}
\end{table*}

\begin{table*}[]
\begin{center}
\begin{tabular}{lcccr}
\toprule
 q\% | \textbf{Expectation} & Test:  (without cv) & cpu\_without & Test: \textbf{CVaR} (with cv) & cpu\_with \\
\midrule
0 & 5.62 & 0.04 & 5.61 & 0.04 \\
25 & 5.63 & 0.04 & 5.62 & 0.05 \\
50 & 5.68 & 0.05 & 5.63 & 0.33 \\
80 & 5.86 & 0.05 & 5.64 & 0.63 \\
90 & 6.03 & 0.05 & 5.75 & 0.63 \\
95 & 6.29 & 0.05 & 5.88 & 0.63 \\
99 & 7.56 & 0.05 & 7.56 & 0.64 \\
100 & 7.64 & 0.05 & 7.64 & 0.65 \\
\bottomrule
\end{tabular}
\e 
\begin{tabular}{lcccr}
\toprule
q\% | \textbf{VaR} & Test: \textbf{CVaR} (without cv) & cpu\_without & Test: \textbf{CVaR} (with cv) & cpu\_with \\
\midrule
0 & 7.82 & 0.04 & 7.81 & 0.04 \\
25 & 8.11 & 0.04 & 8.00 & 0.05 \\
50 & 8.63 & 0.05 & 8.59 & 0.34 \\
80 & 10.12 & 0.05 & 10.07 & 0.63 \\
90 & 11.16 & 0.05 & 11.00 & 0.63 \\
95 & 13.10 & 0.05 & 12.90 & 0.63 \\
99 & 16.24 & 0.05 & 16.24 & 0.64 \\
100 & 32.81 & 0.05 & 32.81 & 0.64 \\
\bottomrule
\end{tabular}
\e 
 \begin{tabular}{lcccr}
\toprule
q\% | \textbf{CVaR} & Test:  (without cv) & cpu\_without & Test: \textbf{CVaR} (with cv) & cpu\_with \\
\midrule
0 & 8.39 & 0.04 & 8.39 & 0.04 \\
25 & 8.84 & 0.04 & 8.83 & 0.04 \\
50 & 9.94 & 0.04 & 9.83 & 0.33 \\
80 & 13.02 & 0.04 & 12.60 & 0.62 \\
90 & 18.02 & 0.04 & 16.68 & 0.62 \\
95 & 24.28 & 0.04 & 24.28 & 0.62 \\
99 & 41.55 & 0.05 & 41.55 & 0.63 \\
100 & 48.94 & 0.05 & 48.94 & 0.63 \\
\bottomrule
\end{tabular}

\end{center}    
\caption{Dispersion statistics (for the different costs) of the optimal solution over the 100 random trials with and without the use of the convex hull-related additional candidate.\\
$\texttt{nDesign}=5$, $N=1000$, $n_\text{cv}=10$, $\texttt{maxIterCV}=10$}\label{resultscvcompar}
\end{table*}

\section{Conclusion}\label{sec-conc}
This paper proposes a unified and efficient algorithm that addresses the problem of stochastic optimization in its different forms. Namely, \texttt{Expectation}, Value-at-risk and Conditional-Value-at-risk-related formulations. 
\e The algorithm leverages a recently proposed \textit{search-\&-accelerate} gradient-based algorithm \cite{alamir2026nonlinearmodelpredictivecontrol}. The latter is used to select rapidly a finite number of candidate solutions which are then compared  via simple enumeration avoiding the problems associated to differentiability and smoothness of some of the formulations. 
\e The trade-of between the quality of the result and the associated computation time might be tailored which might enable the proposed algorithm to be used in a computationally challenging contexts such as the one associated to Stochastic Nonlinear Model Predictive Control (SNMPC) or on-line Stochastic Nonlinear Moving-Horizon Estimation (SNMHE). 
\e The algorithm will be shortly made available for \texttt{python} users through the widely used \texttt{Pypi} repository and can hence be simply accessible through the simple instruction: 
\begin{center}
\texttt{pip install mizoGrad}
\end{center}
which is already the case for its parent algorithm which is dedicated to deterministic gradient-based optimization problems. Moreover, the complete documentation will be shortly available at:
\begin{center}
\texttt{\url{https://www.mazenalamir.fr/mizoGrad/}}
\end{center}
by extending the already existing one that describes the deterministic case. 
\e 
Beside making the algorithm freely available, current undergoing work concerns the incorporation of the algorithm in SNMPC and SNMHE use cases in order to examine its performance in closed-loop. 
\bibliographystyle{IEEEtran}
\bibliography{bib_stochastic.bib}
\end{document}